\documentclass[draft,12pt]{article} 
\usepackage{amsmath,amsfonts,amsthm,amssymb,amscd} 
\renewcommand{\Im}{\mathop{\rm Im}\nolimits} 

\newcommand{\be}{\begin{equation}} 
\newcommand{\ee}{\end{equation}}

\theoremstyle{plain} \newtheorem{theorem}{Theorem}[section] 
\newtheorem{lemma}[theorem]{Lemma} 
\newtheorem{proposition}[theorem]{Proposition} 
 \theoremstyle{definition} 
\newtheorem{definition}[theorem]{Definition} \theoremstyle{remark} 
\newtheorem{remark}[theorem]{Remark}

\newcommand{\R}{{\mathbb R}} \newcommand{\U}{{\mathcal U}}

\newcommand{\Nn}{{\mathbb N}}

\newcommand{\E}{{\mathcal E}} 
\newcommand{\Tr}{{\mathcal T}}

\newcommand{\resto}{{\mathcal 
R}} 
\def\im{{\rm i}} 
\newcommand{\N}{{\mathcal N}}

\newcommand{\V}{{\mathcal V}}

\newcommand{\C}{\mathbb{C}} 
\newcommand{\T}{\mathbb{T}}

\newcommand{\Ze}{{\mathcal Z}} 
\newcommand{\Ye}{{\mathcal Y}} 
\def\norma#1{\left\| #1\right\|}

\def\hsi{H^{\sig}(M,\R)} 
\def\csi{C_{\sig}^\infty(\U,\R)} 
\def\hsic{H^{\sig}(M,\C)} 
\def\lm{\Lambda_m} 
 
\def\uno{\mathbf{1}} 
\def\sig{s} 
 
\def\EL#1#2#3{{\mathcal L}^{#1,#2}_{#3}\null} 
\def\EM#1#2#3{{\mathcal M}^{#1,#2}_{#3}\null} 
\def\ELt#1#2#3{\widetilde{{\mathcal L}}^{#1,#2}_{#3,\ell}\null} 
\def\EMt#1#2#3{\widetilde{{\mathcal M}}^{#1,#2}_{#3,\ell}\null} 
\def\hcs#1{H^{#1}(M,\C)} 

\def\Hs#1{\mathcal{H}^s_{#1}}
 
\numberwithin{equation}{section} 
 
\begin{document} 
 
\author{D. Bambusi, J.-M. Delort, B. Gr{\'e}bert, J. Szeftel}
  
\title{Almost global existence for Hamiltonian semi-linear
Klein-Gordon equations with small Cauchy data on Zoll manifolds}

\date{} 
\maketitle 
 
\begin{abstract} 
This paper is devoted to the proof of almost global existence results
for Klein-Gordon equations on Zoll manifolds (e.g. spheres of
arbitrary dimension) with Hamiltonian
nonlinearities, when the Cauchy data are smooth and small. The proof
relies on Birkhoff normal form methods and on the specific
distribution of eigenvalues of the laplacian perturbed by a potential
on Zoll manifolds.
\end{abstract} 
 
 
\setcounter{section}{0}

\section{Introduction} 
 
Let $(M,g)$ be a compact Riemannian manifold without boundary, denote 
by $\Delta_{g}$ its Laplace-Beltrami operator, and consider the 
nonlinear Klein-Gordon equation  
\be \label{KG} 
(\partial_{t}^2-\Delta_{g}+V+m^2)v=-\partial_{2}f(x,v)  
\ee  
where $m$ is a strictly positive constant, $V$ is a smooth nonnegative 
potential on $M$ and $f\in C^\infty(M\times \R)$ vanishes at least at 
order 3 in $v$, $\partial_{2}f$ being the derivative with respect to 
the second variable. In this work we prove that, for a special class 
of manifolds and for almost every value of $m>0$, this 
\textit{Hamiltonian} partial differential equation admits a Birkhoff 
normal form at any order.  The principal dynamical consequence is the 
almost global existence of small amplitude solutions for such a  
nonlinear Klein-Gordon equation. 
 
More precisely, if $M$ is a Zoll manifold (i.e. a compact manifold 
whose geodesic flow is periodic, e.g. a sphere), for almost every 
value of $m>0$ and for any $N\in \Nn$, we prove that there is $s\gg 1$ 
such that, if the initial data 
$(v\arrowvert_{t=0},\partial_{t}v\arrowvert_{t=0})$ are of size $\epsilon 
\ll 1$ in $H^s \times H^{s-1}$, \eqref{KG} has a solution defined on a 
time interval of length $C_{N}\ \epsilon^{-N}$.  As far as we know, 
this is the first result of that type when the dimension of the 
manifold is larger or equal to 2. 
 
\medskip 
 
Let us recall some known results for the similar problem on $\R^d$, 
when the Cauchy data are smooth, compactly supported, of size 
$\epsilon\ll 1$. In this case, linear solutions decay in $L^\infty$ 
like $t^{-d/2}$ when $t\to\infty$. This allows one to get global 
solutions including quasi-linear versions of (\ref{KG}), when 
$d\geq2$ (see Klainerman \cite{K1} and Shatah \cite{Sh} if $d\geq3$ 
and Ozawa, Tsutaya and Tsutsumi \cite{OTT} if $d=2$). When $d=1$  
Moriyama,  Tonegawa and Tsutsumi \cite{MTT} proved that 
solutions exist over intervals of time of exponential length 
$e^{c/\epsilon^2}$. This result is in general optimal (see references 
in \cite{D}), but global existence for small $\epsilon>0$ was proved 
in \cite{D} when the nonlinearity satisfies a special condition (a ``null 
condition'' in the terminology introduced by Klainerman in the case 
of the wave equation in 3--space dimensions \cite{K}).  
 
For the problem we are studying here, since we have no dispersion on a
compact manifold, we cannot hope to exploit any time decay of the
solutions of the linear equation. Instead we shall use a normal form
method. Remark that if in (\ref{KG}) the nonlinearity vanishes at
order $p\geq2$ at $v=0$, local existence theory gives a solution
defined on an interval of length $c\epsilon^{-p+1}$.  Recently, in
\cite{DS1}, \cite{DS2} Delort and Szeftel proved that the solution of
the same equation exists, for almost all $m>0$, over a time interval
of length $c\epsilon^{-q+1}$, where $q$ is an explicit number strictly
larger than $p$ (typically $q=2p-1$). Actually these papers concern
more general nonlinearities than the one in \eqref{KG}, namely a
suitable class of non Hamiltonian nonlinearities depending on time and
space derivatives of $v$.
 
One of the ideas developed by Delort-Szeftel consists in reducing, by 
normal form procedure, \eqref{KG} to a new system in which the 
nonlinearity vanishes at order $q>p$ at the origin.  In \cite{DS2} an 
explicit computation showed that the first order normal form (which 
leads to a nonlinearity of degree $q$) conserves also the $H^s$ norm 
for any large $s$, whence the result cited above. 
 
\medskip 
 
On the other hand in \cite{BG} Bambusi and Gr{\'e}bert proved an 
abstract Birkhoff normal form theorem for Hamiltonian PDEs. Although 
that theorem remains valid in all dimensions, it supposes that the 
nonlinearity satisfies a ``tame modulus" property. In \cite{BG} this 
property was only verified for a quite general class of $1-d$ PDEs and 
for a particular NLS equation on the torus $\T^d$ with arbitrary $d$. 
Actually in that paper, the tame modulus property was verified by the 
use of the property of ``well localization with respect to the 
exponentials" established by Craig and Wayne \cite{CW}, a property which 
has no equivalent in higher dimensions. 
 
It turns out that in \cite{DS2} Delort and Szeftel proved an estimate  
concerning multilinear forms defined on $M$ that implies a weaker form 
of the tame modulus property assumed in \cite{BG}.  
The present paper is the result of the combination of the arguments  
of \cite{DS1}, \cite{DS2} and of \cite{BG}.

\medskip 
 
We recall that some other partial normal form results for PDEs have been 
previously obtained by Kuksin and P{\"o}schel \cite{KP96}, by Bourgain 
\cite{Bo96,Bo04} and, for perturbations of completely resonant systems, by 
Bambusi and Nekhoroshev \cite{BN98}.  For a more precise discussion we 
refer to the introduction of \cite{BG}.

\medskip 
 
Let us conclude this introduction mentioning several open questions.
The first concerns the possibility of proving almost global existence
for more general nonlinearities than the Hamiltonian ones we consider
here. Of course, one cannot expect to be able to do so for any
nonlinearity depending on $v$ and its first order derivatives: in
\cite{D1} an example is given on the circle ${\mathbb{S}}^1$ of a
nonlinearity for which the solution does not exist over a time
interval of length larger than the one given by local existence theory
(Remark that this example holds true for any value of $m>0$). On the
other hand, Delort and Szeftel constructed in \cite{DS3} almost global
solutions of equations of type \eqref{KG} on manifolds of revolution,
for radial data, with a nonlinearity $f$ depending on
$(v,\partial_{t}v)$ and even in $\partial_{t}v$. We thus ask the
question of finding a ``null condition" (in the spirit of Klainerman
\cite{K}) for semi-linear nonlinearities $f(v,\partial_{t}v, \nabla
v)$, which would allow almost global existence of small $H^s$
solutions for almost every $m>0$.
 
The second question we would like to mention concerns the  
exceptional values of $m$ which are excluded of our result.  The  
conservation of the Hamiltonian of equation \eqref{KG} allows one to  
control the $H^1$-norm of small solutions. This implies global  
existence of small $H^1$ solutions in one or two space dimensions. The  
results we establish in the present paper show that for almost every  
$m>0$, the $H^s$-norms of these solutions remain small over long time  
intervals if they are so at $t=0$. What happens when $m$ is in the  
exceptional set? In \cite{Bo96b} Bourgain constructed, in one space  
dimension and for a convenient perturbation of $-\Delta$, an example  
of a solution whose $H^s$-norm grows with time. Nothing seems to be  
known in larger dimensions. In particular, if $d\geq 3$, one does not  
even know if for all $m>0$ a solution exists almost globally,  
eventually without staying small in $H^s$ ($s\gg 1$).

\section{Statement of main results} 
 
We begin, in section \ref{subsec1.1}, by a precise exposition of our 
result concerning the almost globality. The Birkhoff normal form 
theorem for equation \eqref{KG} that implies the almost globality 
result will be presented in  section \ref{subsec1.3}, after the 
introduction of the Hamiltonian formalism in section \ref{subsec1.2}. 
 
\subsection{Almost global solution}\label{subsec1.1} 

Let $(M,g)$ be a compact Riemannian manifold without boundary of
dimension $d\geq 1$. Denote by $\Delta_{g}$ its Laplace-Beltrami
operator.  Let $V$ be a smooth nonnegative potential on $M$ and $m\in
(0,\infty)$. Let $f\in C^\infty(M\times \R)$ be such that $f$ vanishes at
least at order 3 in $v$. We consider the following Cauchy problem for
the nonlinear Klein-Gordon equation
\begin{align} 
\begin{split} \label{KGCauchy} 
(\partial_{t}^2-\Delta_{g}+V+m^2)v&=-\partial_{2}f(x,v) 
\\ 
v\arrowvert_{t=0}&=\epsilon v_{0} 
\\ 
\partial_{t}v\arrowvert_{t=0}&=\epsilon v_{1} 
\end{split} 
      \end{align} 
where $v_{0}\in H^s(M,\R)$, $v_{1}\in H^{s-1}(M,\R)$ are real valued  
given data and $\epsilon>0$. We shall prove  
that the above problem has almost global solutions for almost every  
$m$ when $\epsilon>0$ is small enough and $s$ is large enough, under  
the following geometric assumption on $M$: 
\definition One says that $(M,g)$ is a Zoll manifold if and only if  
the geodesic flow is periodic on the cosphere bundle of $M$. 

\medskip
 
Our main dynamical result is the following: 
 
\begin{theorem} \label{thm1} 
Let $(M,g)$ be a Zoll manifold and let $V:M\to \R$ be a smooth 
nonnegative potential. Let $r\in \Nn$ be an arbitrary integer. There 
is a zero measure subset $\mathcal N$ of $(0,+\infty)$, and for any  
$m\in (0,+\infty) \setminus \mathcal N$, there is $s_{0}\in\Nn$ such that  
for any $s\geq s_{0}$, for any real valued $f\in C^\infty(M\times \R)$ 
vanishing at least at order 3 at $v=0$, there are $\epsilon_{0}>0$, 
$c>0$, such that for any pair $(v_{0},v_{1})$ of real valued functions 
belonging to the unit ball of $H^s(M,\R)\times H^{s-1}(M, \R)$, any 
$\epsilon \in (0,\epsilon_{0})$, the Cauchy problem \eqref{KGCauchy} 
has a unique solution 
$$ v\in C^0((-T_{\epsilon},T_{\epsilon}),H^s(M,\R))\cap
C^1((-T_{\epsilon},T_{\epsilon}),H^{s-1}(M,\R))$$ with
$T_{\epsilon}\geq c\epsilon^{-r}$. Moreover there is $C>0$ such that,
for any $t\in (-T_{\epsilon},T_{\epsilon})$, one has
\be
\label{estimHs} \Vert{v(t,\cdot )}\Vert_{H^s}+\Vert\partial_{t}
v(t,\cdot )\Vert_{H^{s-1}}\leq C\epsilon \ .  
\ee
\end{theorem} 
 
\noindent 
{\bf Comments} 
 
The above theorem provides Sobolev bounded almost global solutions for 
equation \eqref{KGCauchy} with small smooth Cauchy data on a 
convenient class of compact manifolds. To our knowledge this is the 
first result of this kind on compact manifolds of dimension larger or equal  
to 2. In the case of one dimensional compact manifolds, similar 
statements have been obtained by Bourgain \cite{Bo96,Bo04} (with a loss on 
the number of derivatives of the solution with respect to those of the 
data), by Bambusi \cite{Bam03} and by Bambusi-Gr{\'e}bert 
\cite{BG}. Remark that in this case, because of the conservation of 
the Hamiltonian of the equation, one controls uniformly the $H^1$-norm 
of small solutions, which implies global existence of such 
solutions. The results of the preceding authors allow to control 
$H^s$-norms of these solutions for very long times. In the case of  
compact manifolds of revolution and for convenient radial data,  
Delort and Szeftel got in \cite{DS3} Sobolev bounded almost global 
solutions (remark that this result is morally one-dimensional). 
 
\medskip 
 
The assumption that $M$ is a Zoll manifold will be used in the proof  
through distribution properties of the eigenvalues of the Laplacian  
of $M$. Actually we shall prove theorem \ref{thm1} for any  
compact manifold without boundary ($M,g)$ such that if 
\be \label{P} 
P=\sqrt{-\Delta_{g} +V}, 
\ee 
the spectrum $\sigma(P)$ of $P$ satisfies the following condition:  
there are constants $\tau >0$, $\alpha \in \R$, $c_0 >0$, $\delta  
>0$, $C_0 >0$, $D\geq 0$, and a family of disjoint compact  
intervals $(K_{n})_{n\geq 1}$, with $K_{1}$ at the left of $K_{2}$  
and for $n\geq 2$ 
\be \label{Kn} 
K_{n}=\left[ \frac{2\pi}{\tau}n+\alpha -\frac{c_{0}}{n^\delta}, 
\frac{2\pi}{\tau}n+\alpha +\frac{c_{0}}{n^\delta}\right], 
\ee 
such that  
\begin{align}\begin{split}\label{115} 
\sigma(P)&\subset \bigcup_{n\geq 1}K_{n}\\ 
\#(\sigma(P)\cap K_{n})&\leq C_{0}n^D\ . 
\end{split}\end{align} 
If $M$ is a Zoll manifold, and if $\tau >0$ is the minimal period of  
the geodesic flow on $M$, the results of Colin de Verdi{\`e}re  
\cite{CV} (see also Guillemin \cite{G} and Weinstein \cite{W}) show  
that the large eigenvalues of $P$ are contained inside the union of  
the intervals 
$$ 
\left[ \frac{2\pi}{\tau}n+\alpha -\frac{C}{n}, 
\frac{2\pi}{\tau}n+\alpha +\frac{C}{n}\right] 
$$ 
for $n$ large enough and for some constant $C>0$. Making a translation  
in $n$ and $\alpha$, and changing the definition of the constants, one  
sees that this implies conditions \eqref{Kn}, \eqref{115} for any  
$\delta \in (0,1)$ (remark that the second condition in \eqref{115}  
holds true with $D=d-1$ because of Weyl law). 
 
On the other hand conditions \eqref{Kn}, \eqref{115} are not more  
general than the assumption that $M$ is a Zoll manifold, since by  
theorem 3.2 in Duistermaat and Guillemin~ \cite{DG}, they imply  
that the geodesic flow is periodic. 
 
\subsection{Hamiltonian formalism}\label{subsec1.2} 
 
We introduce here (see e.g. \cite{CheMa}) the Hamiltonian 
formalism we shall use to solve the equation. We denote by 
\begin{equation} 
\label{1.2.1} 
\langle f_1,f_2\rangle 
\end{equation} 
the bilinear pairing between complex valued distributions and test 
functions on $M$. We shall use the same notation for vector valued 
$f_1,f_2$.  
 
If $F$ is a $C^\infty$ function on an open subset $\U$ of the Sobolev 
space of real valued functions $\hsi$, $\sig\geq 0$, we define for 
$p\in\U$, the $L^2$ gradient $\nabla F(p)$ by 
\begin{equation} 
\label{1.2.2} 
\partial F(p)h= \langle\nabla F(p), h \rangle\ ,\quad \forall h \in \hsi, 
\end{equation} 
 $\partial F$ denoting the differential. In that way $\nabla F(p) $ is 
an element of $H^{-\sig}(M,\R)$. When we consider real valued 
$C^\infty$ functions defined on an open subset of 
$\hsi\times \hsi\equiv \hsi^2$, $(p,q)\mapsto F(p,q)$ we write 
\begin{eqnarray*} 
\partial F(p,q)&=& (\partial_p F(p,q),\partial_q F(p,q) ) 
\\ 
\nabla F(p,q)&=& (\nabla_p F(p,q),\nabla_q F(p,q) )\in 
H^{-\sig}(M,\R)\times H^{-\sig}(M,\R). 
\end{eqnarray*} 
 
Endow $\hsi^2$ with the weak symplectic structure 
\begin{equation} 
\label{1.2.4} 
\Omega\left((p,q),(p',q')\right):= \langle q,p' 
\rangle - \langle q',p\rangle  = \langle J^{-1}(p,q),(p',q') \rangle 
\end{equation} 
where $J$ is given by 
\begin{equation} 
\label{1.2.5} 
J= 
\left[ 
\begin{matrix} 
0 & -\uno 
\\ 
\uno &0 
\end{matrix} 
\right]. 
\end{equation} 
If $\U$ is an open subset of $\hsi^2$ and $F\in C^\infty (\U,\R)$, then, 
for $(p,q)\in\U$, we define its Hamiltonian vector field by 
\begin{equation} 
\label{1.2.6} 
X_F(p,q)=J\nabla F(p,q)=(-\nabla_qF(p,q),\nabla_pF(p,q)) 
\end{equation} 
which is characterized by 
\begin{equation} 
\label{1.2.7} 
\Omega\left( X_F,(h_p,h_q)\right)=\partial F (h_p,h_q)= \partial_p F 
h_p+\partial_q F h_q  
\end{equation} 
for any $(h_p,h_q)\in\hsi^2$. 
 
A special role is played by the functions whose Hamiltonian vector 
field is an $\hsi^2$ valued function. Thus we give the following 
 
\begin{definition} 
\label{d.1.2.1} 
If $\U$ is an open subset of $\hsi^2$, we denote by $\csi$
(resp. $C^\infty_\sig(\U,\C)$) the space of real (resp. complex)
valued $C^\infty$ functions defined on $\U$ such that
\begin{equation} 
\label{1.2.8} 
X_F\in C^\infty(\U,\hsi^2 )\quad (\text{or}\quad \nabla F\in 
C^\infty(\U,\hsi^2 ) ), 
\end{equation}  
resp. 
\begin{equation} 
\label{1.2.8bis} 
X_F\in C^\infty(\U,\hsi^2 \otimes\C)\quad (\text{or}\quad \nabla F\in 
C^\infty(\U,\hsi^2 \otimes\C) ). 
\end{equation}  
\end{definition} 
We shall use complex coordinates in $\hsi^2$ identifying this space 
with $\hsic$, through $(p,q)\mapsto u=(p+\im q)/\sqrt2$. We set 
\begin{eqnarray} 
\label{1.2.9} 
&\partial_u=\frac{1}{\sqrt2}(\partial_p-\im\partial_q),
&\partial_{\bar u}=\frac{1}{\sqrt2}(\partial_p+\im\partial_q) 
\\ 
\label{1.2.9.1} 
&\nabla_u=\frac{1}{\sqrt2}(\nabla_p-\im\nabla_q),
&\nabla_{\bar u}=\frac{1}{\sqrt2}(\nabla_p+\im\nabla_q) 
\end{eqnarray}   
so that, if $F$ is a $C^1$ real valued function, we have an identification 
\begin{equation} 
\label{1.2.10} 
X_F(u,\bar u)=\im \nabla_{\bar u}F(u,\bar u)\ . 
\end{equation} 
 
If $F\in\csi$, then clearly $X_F\in C^\infty(\U,\hsic)$.

 For $m \in (0,+\infty)$ let us define 
\begin{equation}\label{1.2.3} 
\Lambda_m = \sqrt{-\Delta_g+V+m^2}. 
\end{equation} 
Let $\sig>(d-1)/2$. We shall write equation \eqref{KG} as a 
Hamiltonian system for $p=\lm^{-1/2}\partial_t v $ and $q=\lm^{1/2}v$ 
on $\hsi^2$. Define 
\begin{equation} 
\label{1.2.11} 
G_2(p,q)=\frac{1}{2}\int_M 
\bigl(\bigl|\lm^{1/2}p\bigr|^2+\bigl|\lm^{1/2}q\bigr|^2\bigr) dx\ 
,\quad \tilde G(p,q)=\int_M f(x,\lm^{-1/2}q) dx 
\end{equation} 
where $dx$ is the Riemannian volume on $M$, and set 
\begin{equation} 
\label{1.2.12} 
G=G_2+\tilde G. 
\end{equation} 
Then by (\ref{1.2.6}) 
\begin{eqnarray} 
\label{1.2.13} 
X_{G_2}(p,q)=(-\lm q,\lm p)\ ,\quad X_{\tilde G}(p,q)= (-\lm^{-1/2} 
\partial_2f(x,\lm^{-1/2} q),0) 
\end{eqnarray} 
where $\partial_2f$ is the derivative with respect to the second 
argument. Then one 
has that $\tilde G\in \csi$ with $\U=\hsi^2$ (actually $X_{\tilde G}$ takes 
values in $H^{\sig+1}(M,\R)^2$).  
 
It follows also that equation (\ref{KG}) can be written as 
\begin{equation} 
\label{1.2.15} 
(\dot p,\dot q)=X_G(p,q) 
\end{equation} 
or, using (\ref{1.2.10}) 
\begin{equation} 
\label{1.2.16} 
\dot u=\im \nabla_{\bar u}G(u,\bar u). 
\end{equation} 
 
In the rest of this section we shall give a few technical results that 
we shall need for the proofs of theorems \ref{thm1}, \ref{thm2}. 
 
\begin{definition} 
\label{d.1.2.2} 
Let $\U$ be an open subset of $\hsi^2$ and $F_j\in\csi$, 
$j=1,2$. Then their Poisson bracket is defined by 
\begin{equation} 
\label{1.2.17} 
\left\{F_1,F_2\right\}=\partial F_2 
\cdot X_{F_1}=\Omega(X_{F_2},X_{F_1}) 
\end{equation} 
and one has $\{F_1,F_2\}\in\csi $. 
 
One extends the definition to complex valued functions by linearity of 
the bracket relatively of each of its arguments. 
\end{definition} 
 
The fact that (\ref{1.2.17}) has a smooth vector field follows from 
the well known formula 
\begin{equation} 
\label{1.2.17.1} 
X_{\{ F_1,F_2\}}=[X_{F_1},X_{F_2}]=\partial X_{F_2}\cdot X_{F_1}-\partial 
X_{F_1}\cdot X_{F_2}\ , 
\end{equation} 
with the square bracket denoting the Lie bracket of vector fields (for 
a proof of this formula in the case of weak symplectic manifolds see 
\cite{Bam93}). In case either $F_1$ or $F_2$ do not have a smooth 
vector field, one can also define their Poisson brackets by formula 
(\ref{1.2.17}) but one has to check that it is a well defined 
function, using the fact that we may write 
\begin{equation} 
\begin{split} 
\label{1.2.18} 
\left\{F_1,F_2\right\} &= - (\partial_p F_2)(\nabla_q F_1) +
(\partial_q F_2)(\nabla_p F_1)\\ &= - \langle\nabla_p F_2,\nabla_q F_1
\rangle + \langle\nabla_q F_2,\nabla_p F_1 \rangle \\ &= \im
(\partial_u F_2)(\nabla_{\bar u} F_1) -\im (\partial_{\bar u}
F_2)(\nabla_u F_1).
\end{split} 
\end{equation} 
Let us recall also the rule of transformation of vector fields and 
Poisson brackets under  symplectomorphism. Let $\U$ and $\V$ be open 
subsets of $\hsi^2$, and $\chi:\U\to\V$ be a smooth symplectic 
diffeomorphism. We have by definition for any $u \in \U$ 
\begin{equation}\label{1.2.19} 
(\partial\chi(u))^{-1} = J{}^t(\partial\chi(u))J^{-1}. 
\end{equation} 
For $F\in C^{\infty}_\sig(\V,\R)$ one has 
\begin{equation} 
\label{1.2.20} 
X_{F\circ \chi}(u)=(\partial \chi(u))^{-1}X_F(\chi(u)) 
\end{equation} 
and therefore $F\circ\chi\in C^{\infty}_\sig(\U,\R) $ (actually
(\ref{1.2.20}) holds in the more general context where $\nabla F$ has
a domain which is left invariant by $\chi$). We also remark that for
any $C^1$ real-valued function $F_1$ on $\V$ and for any $F_2$ in
$C^{\infty}_\sig(\V,\R)$ one has
\begin{equation} 
\label{1.2.21} 
\left\{F_1\circ\chi,F_2\circ\chi\right\}= 
\left\{F_1,F_2\right\}\circ\chi . 
\end{equation}

To conclude this subsection let us state as a lemma the well known 
formula that is the root of the Birkhoff normal form method as 
developed using Lie transform.
\begin{lemma} 
\label{l.1.2.3} 
Let $F,G$ be two real valued functions defined on 
$\U\subset\hsi^2$. Assume that $F\in C^\infty_\sig(\U,\R)$ and 
$G\in C^\infty(\U,\R)$. Denote by 
$(Ad\, F)\, h=\left\{F,h\right\}$. Then $(Ad\,F)G$ is well
defined, and if we assume that  for some $n\geq 1$
\begin{equation} 
\label{1.2.22} 
F_n:= (Ad\,F)^n G
\end{equation}   
is well defined and belongs to  $C^\infty_\sig(\U,\R)$, then $F_{n+1}$
is also well defined.

Let $\V$ be such that $\overline{\V}\subset \U$. There 
exists a positive $T$ such that the flow $\V\ni(p,q)\mapsto 
\Phi^t(p,q)\in\U$ of $X_F$ is well defined and smooth for 
$|t|<T$. Moreover, for $|t|<T$ and $(p,q)\in\V$, one has for any $r\in
\Nn$ the formula 
\begin{equation} 
\label{1.2.23} 
G(\Phi^t(p,q))=\sum_{n=0}^{r}\frac{t^n}{n!}F_n(p,q)+\frac{1}{r!}\int 
_0^t (t-s)^{r}F_{r+1}(\Phi^s(p,q))ds. 
\end{equation} 
\end{lemma} 
\proof Remark first that $(Ad\,F)G$ is well defined by \eqref{1.2.18},
and that under our assumptions, for $n\geq 2$, $F_n$ is well defined
by definition~\ref{d.1.2.2}. Since $X_F$ is smooth on $\U$ the flow
$\Phi^t(.)$ is a smooth symplectic diffeomorphism on $\V$. For fixed
$(p,q)$ put $\phi(t)=G(\Phi^t(p,q))$. Formula (\ref{1.2.23}) follows
from Taylor formula since $\phi(t)$ is $C^\infty$. We thus have
$\phi'(t)=[(Ad\, F)\, G ](\Phi^t(p,q))=F_1(\Phi^t(p,q))$. Using
(\ref{1.2.22}) one proves by induction that
$\phi^{(n)}(t)=F_n(\Phi^t(p,q))$ and the conclusion follows. \qed

\subsection{Birkhoff Normal Form} 
\label{subsec1.3} 
 
Using the notation of section \ref{subsec1.1}, we define for $n\geq 1$  
spectral projectors 
\be \label{pin} 
\Pi_{n}=\mathbf{1}_{K_{n}}(P)\ . 
\ee 
Then, for $(p,q)\in H^{s}(M,\R)^2$  we introduce the quantities 
\be \label{Jn} 
J_{n}(p,q)=\frac{1}{2}\left( \Vert\Pi_{n}p\Vert^2_{L^2}  
+\Vert\Pi_{n}q\Vert^2_{L^2}\right)\ . 
\ee 

For $(p,q)\in H^{s}(M,\R)^2$ we denote 
$$ 
\norma{(p,q)}_{s}^2:= \Vert p\Vert^2_{H^s}+\Vert q\Vert^2_{H^s}\ 
$$ 
 
We can now state our Birkhoff normal form result for the nonlinear  
Klein-Gordon equation on Zoll manifolds:

\begin{theorem}\label{thm2} 
    Let $G$ be the Hamiltonian given by \eqref{1.2.11},
    \eqref{1.2.12}.  Then for any $r\geq 1$, there exists a zero
    measure subset $\mathcal N$ of $(0,+\infty)$, and for any $m\in
    (0,+\infty) \setminus \mathcal N$, there exists a large $s_0$ with
    the following properties: For any $s\geq s_0$, there exist two
    neighborhoods of the origin $\U$, $\V$, and a bijective canonical
    transformation $\Tr:\V\to\U$ 
    which puts the Hamiltonian in the form \begin{equation}
    \label{for.nor} G\circ\Tr= G_{2}+\Ze+\resto \end{equation} where
    $\Ze$ is a real valued continuous polynomial of degree at most
    $r+2$ satisfying
\begin{equation} 
\label{stime1} 
\left\{J_n,\Ze \right\}=0\ ,\quad \forall n\geq 1 
\end{equation} 
and $\resto\in C^\infty_s(\V,\R)$ has a zero of order $r+3$ at 
the origin. Precisely its vector field fulfills the estimate 
    \begin{equation} 
    \label{stime2} 
    \norma{X_{\resto}(p,q)}_s\leq 
    {C_s} \norma{(p,q)}_s^{r+2}\ ,\quad (p,q)\in\V.
    \end{equation} 
    Finally the canonical transformation satisfies
    \begin{equation} 
    \label{stime} 
    \norma{(p,q)-\Tr(p,q)}_s\leq C_s 
    \norma{(p,q)}_s^2\ ,\quad (p,q)\in\V.
    \end{equation} 
    Exactly the same estimate is fulfilled on $\U$ by the inverse canonical 
    transformation. 
\end{theorem} 

From \eqref{stime} it follows $\Tr(0)=0$ and $\partial\Tr(0)=\uno $.

Theorem \ref{thm2} implies theorem \ref{thm1} (see the proof of theorem 
\ref{thm1} in section \ref{sec:proofthm1}) but it says more: 
namely, the $J_{n}$ are almost conserved  
quantities for the equation \eqref{KG}. More 
precisely, with the notation of theorems \ref{thm1} and \ref{thm2}, 
for any $n\geq 1$  
\be \label{estimJn} \vert J_{n}(p(t),q(t))- 
J_{n}(p(0),q(0))\vert\leq \frac C{n^{2s}}\epsilon^3 
\quad \mbox{ for } |t|\leq \epsilon^{-r}  
\ee  
where $p(t)=\lm^{-1/2}\partial_t v(t) $ and $q(t)=\lm^{1/2}v(t)$ (for
the proof see the end of section \ref{sec:proofthm1}). 
Roughly speaking, the last property means that energy transfers are 
allowed only between modes corresponding to frequencies in the same 
spectral interval $K_{n}$.

\section{Proof of the main results} 
\label{proof} 
In this section we prove theorem \ref{thm2} and then deduce theorem 
\ref{thm1}. The proof uses a Birkhoff procedure described in subsection 
\ref{birk}. 
Formally this procedure is very close to the classical Birkhoff scheme
in finite dimension. Nevertheless, in infinite dimension, we need to
define a convenient framework in order to justify the formal
constructions.  This framework, first introduced in \cite{DS2}, is
presented, and adapted to our context, in the next subsection.

\subsection{Multilinear Forms} 
\label{multi} 
 
Let us introduce some notations. If $n_1,\ldots,n_{k+1}$ are in $\Nn^*$, 
we denote the second and third largest elements of this family by 
\begin{equation} 
\begin{split}
\label{2.1.1} 
\max\!{}_2(n_1,\ldots,n_{k+1})&=\max\left(\left\{n_1,\ldots,n_{k+1}\right\} 
-\left\{n_{i_0}  
\right\}\right) 
\\ 
\mu (n_1,\ldots,n_{k+1})&=\max\left(\left\{n_1,\ldots,n_{k+1}\right\} 
-\left\{n_{i_0},n_{i_1}  
\right\}\right) 
\end{split}
\end{equation} 
where $i_0$ and $i_1$ are the indices such that  
$$ n_{i_0}=\max\!{}(n_1,\ldots,n_{k+1})\ ,\quad 
n_{i_1}=\max\!{}_2(n_1,\ldots,n_{k+1})
$$ 
 and where by convention, when $k=1$, 
$\mu (n_1,n_2)=1$.
We define then  
\begin{equation} 
\begin{split}
\label{2.1.2} 
S(n_1,\ldots,n_{k+1})&=\sum_{\ell=1}^{k+1}[n_\ell-\sum_{j\not=\ell}n_j]_++ 
\mu (n_1,\ldots,n_{k+1}) 
\end{split}
\end{equation} 
where $[a]_+=\max(a,0)$. If $n_k$ and $n_{k+1}$ are the largest two among 
$n_1,\ldots,n_{k+1}$, we have 
\begin{equation} 
\begin{split}
\label{2.1.3} 
\mu (n_1,\ldots,n_{k+1})&\sim n_1+\cdots+n_{k-1} +1
\\  
S(n_1,\ldots,n_{k+1})&\sim \left|n_k-n_{k+1}\right|+n_1+\cdots+n_{k-1}+1\ . 
\end{split}
\end{equation} 
We shall denote by $\E$ the algebraic direct sum of the ranges of the 
$\Pi_n$'s defined by (\ref{pin}). 
 
\begin{definition} 
\label{d.2.1.1} 
Let $k\in\Nn^*$, $\nu\in[0,+\infty)$, $N\in\Nn$.
\begin{itemize} 
\item[i)] We denote by $\EL\nu N{k+1}$ the space of $(k+1)$--linear 
  forms $L:\E\times \cdots\times \E\to\C$ for which there exists $C>0$ such 
  that for any $u_1,\ldots,u_{k+1}\in\E$, any $n_1,\ldots,n_{k+1}$ in $\Nn^*$ 
\begin{equation} 
\label{2.1.4} 
\left|L(\Pi_{n_1}u_1,\ldots,\Pi_{n_{k+1}}u_{k+1})  \right|\leq 
C\frac{\mu (n_1,\ldots,n_{k+1})^{\nu+N}} 
{S(n_1,\ldots,n_{k+1})^N}\prod_{j=1}^{k+1} \norma{u_j}_{L^2}. 
\end{equation} 
\item[ii)] We denote by $\EM\nu N{k}$ the space of $k$--linear maps
$M:\E\times \cdots \times \E\to L^2(M,\C)$ for which there exists
$C>0$ such that for any $u_1,\ldots,u_{k}\in\E$ any
$n_1,\ldots,n_{k+1}$ in $\Nn^*$
\begin{equation} 
\label{2.1.5} 
\left\Vert \Pi_{n_{k+1}} M(\Pi_{n_1}u_1,\ldots,\Pi_{n_{k}}u_{k})  
\right\Vert_{L^2} 
\leq C\frac{\mu (n_1,\ldots,n_{k+1})^{\nu+N}} 
{S(n_1,\ldots,n_{k+1})^N}\prod_{j=1}^{k} \norma{u_j}_{L^2}. 
\end{equation} 
\end{itemize} 
The best constant $C$ in (\ref{2.1.4}), (\ref{2.1.5}) defines a norm 
on the above spaces. We set also $\EL\nu {+\infty}{k+1} =
\bigcap_{N\in\Nn} \EL\nu N{k+1}$.
\end{definition} 
 
Consider $L\in \EL \nu N{k+1}$ with $N>1$ and fix an integer 
$j\in\{1,\ldots,k+1\}$ and elements $u_\ell\in\E$ for 
$\ell\in\left\{1,\ldots,k+1\right\} -\left\{j\right\}$. Then by 
\eqref{2.1.3}, \eqref{2.1.4}  
$$ 
\sum_{n_j}L(u_1,\ldots,u_{j-1},\Pi_{n_j}u_j, 
u_{j+1},\ldots,u_{k+1}) 
$$  
converges for any $u_j\in L^2(M,\C)$, so $u_j\mapsto 
L(u_1,\ldots,u_{k+1})$ extends as a continuous linear form on 
$L^2(M,\C)$. Consequently, there is a unique element 
$M_{L,j}(u_1,\ldots,\widehat {u_j},\ldots,u_{k+1})$ of $L^2(M,\C)$ with 
\begin{equation} 
\label{2.1.6} 
L(u_1,\ldots,u_{k+1})=\left\langle u_j, M_{L,j}(u_1,\ldots,\widehat 
  {u_j},\ldots,u_{k+1})  \right\rangle 
\end{equation} 
for all $u_1,\ldots, u_{k+1}\in\E$. By \eqref{2.1.4}, $M_{L,j}$ satisfies 
\eqref{2.1.5}, i.e. defines an element of $\EM\nu Nk$. Conversely, if 
we are given an element of $\EM\nu Nk$, we define a multilinear form 
belonging to $\EL\nu N{k+1}$ by a formula of type \eqref{2.1.6}. 

The 
basic example satisfying definition \ref{d.2.1.1} is provided by the 
following result proved in \cite{DS2} (proposition 1.2.1).
\begin{proposition} 
\label{p.2.1.2} 
Let $k\in \Nn^*$. Denote by $dx$ any measure on $M$ with a $C^\infty$ 
density with respect to the Riemannian volume. There is 
$\nu\in(0,+\infty)$ such that the map 
\begin{equation} 
\label{2.1.7} 
(u_1,\ldots,u_{k+1})\mapsto \int_M u_1\cdots u_{k+1} dx 
\end{equation} 
defines an element of $\EL\nu{+\infty}{k+1}$. 
\end{proposition}

\begin{remark} 
\nonumber Up to now we did not use the spectral assumption 
\eqref{115} on the manifold $M$. Actually proposition 1.2.1 of 
\cite{DS2} is proved on any compact manifold without boundary, 
replacing in \eqref{2.1.4} the spectral projectors $\Pi_n$ defined in 
\eqref{pin} by  spectral projectors $\Pi_\lambda$ associated to 
arbitrary intervals of center $\lambda$ and length $O(1)$. 
\end{remark} 
 
We now use the fundamental example given by the previous proposition
to verify that the nonlinearity $\tilde G$ defined in \eqref{1.2.11}
is in a good class of Hamiltonian functions. If $L$ is a
$(k+1)$-linear map, and if $a\in \Nn$ satisfies $0\leq a\leq k+1$, we
set for $u,\bar u \in \mathcal E$
\be \label{2.1.10} 
\underline L^a(u,\bar u)= L(u,\ldots,u,\bar u,\ldots,\bar u) 
\ee 
where in the right hand side one has $a$ times $u$ and $(k+1-a)$--times  
$\bar u$.  We then define the following  class of Hamiltonian functions:
\begin{definition} 
\label{multi.d} 
For $k \in \Nn$ and $s, \nu \in\R$ with $s>\nu + \frac{3}{2}$, we define 
$\mathcal{H}_{k+1}^s(\nu)$ as the
space of all real valued smooth functions defined on $H^s(M,\C)$,
$(u,\bar{u}) \to Q(u,\bar{u})$, such that there are  for $\ell = 0,\ldots,k+1$
multilinear forms $L_\ell \in \mathcal{L}_{k+1}^{\nu,+\infty}$ with
\[Q(u,\bar{u}) = \sum_{\ell =0}^{k+1} \underline{L}_\ell^\ell(u,\bar{u}).\] 
\end{definition} 
 
This definition is obtain by adapting to our context the usual
definition of polynomial used for example in the theory of analytic
functions on Banach spaces (see for example \cite{Muj} or
\cite{Nik86}).

As a consequence of proposition \ref{p.2.1.2} one gets:
\begin{lemma}\label{lem:G}
Let $P$ be the Taylor's polynomial of $\tilde G$ at degree $k$.
Then there exists 
$\nu\in(0,+\infty)$ such that $P$ can be decomposed as
$$P=\sum_{j=3}^{k} P_{j}$$
where $P_{j} \in \mathcal{H}_{j}^s(\nu)$.
\end{lemma}

Let us recall the main properties for $\EM\nu N{k}$ established  
in proposition 2.1.3 and theorem 2.1.4 of \cite{DS2}. 
\begin{proposition} 
\label{p.2.1.3.bis} 
\begin{itemize} 
\item[i)] Let $\nu\in[0,+\infty)$, $s\in\R$, $s>\nu+3/2$, $N\in\Nn$, 
  $N>s+1$. Then, any element 
  $M\in\EM\nu N{k}$ extends as a bounded operator from 
  ${\hcs s}^{k}$ to $\hcs s$.  
  Moreover, for any $s_0\in(\nu+3/2,s]$, there is $C>0$ such that for any 
  $u_\ell\in\hcs s$, $\ell\in\{1,\ldots,k\}$  
\begin{equation} 
\label{2.1.8.bis} 
\norma{M(u_1,\ldots,u_{k})}_{H^s}\leq C \norma{M}_{\EM\nu
  N{k}}\Big(\sum_{{1\leq \ell\leq
  k}}\norma{u_\ell}_{H^s}\prod_{{\ell'\not=\ell}}
  \norma{u_{\ell'}}_{H^{s_0}}\Big).
\end{equation} 
\item[ii)] Let $k_1,k_2\in\Nn^*$, $\nu_1,\nu_2\in[0,+\infty)$, $1\leq 
  \ell\leq k_2$. For $M_1\in\EM{\nu_1}N{k_1}$, $M_2\in\EM{\nu_2}N{k_2}$ 
  with $N> 1+\max(\nu_1,\nu_2)$, define a $(k_1+k_2-1)$--linear operator on 
  $\E^{k_1+k_2-1}$  
$$ 
(u_1,\ldots,u_{k_1+k_2-1})\to M(u_1,\ldots,u_{k_1+k_2-1}) 
$$  
by 
\begin{equation} 
\label{2.1.9.bis} 
\begin{split}
M(u_1,\ldots,u_{k_1+k_2-1})= \makebox[7cm]{}\\
M_2(u_1,\ldots,u_{\ell-1},M_1(u_\ell,\ldots,u_{\ell+k_1-1}),
u_{\ell+k_1},\ldots,u_{k_1+k_2-1})\ . 
\end{split}
\end{equation} 
Then $M$ belongs to $\EM{\nu_1+\nu_2+1}{N-\max(\nu_1,\nu_2)-1}{k_1+k_2-1}$ 
and the map $(M_1,M_2)\mapsto M$ is bounded from 
$\EM{\nu_1}N{k_1}\times \EM{\nu_2}N{k_2}$ to the preceding space. 
\end{itemize} 
\end{proposition} 
 
Using the duality formula \eqref{2.1.6}, proposition \ref{p.2.1.3.bis}  
immediately implies the corresponding properties for the multilinear forms  
of $\EL\nu N{k+1}$.  
 
\begin{proposition} 
\label{p.2.1.3} 
\begin{itemize} 
\item[i)] Let $\nu\in[0,+\infty)$, $s\in\R$, $s>\nu+3/2$, $N\in\Nn$,
  $N>s+1$. Then for any $j\in\{1,\ldots,k+1\}$, any multilinear form
  $L\in\EL\nu N{k+1}$ extends as a continuous multilinear form
  $(u_1,\ldots,u_{j},\ldots,u_{k+1})\mapsto L
  (u_1,\ldots,u_{j},\ldots,u_{k+1}) $ on
$$ \hcs s\times \cdots\times \hcs s\times \hcs{-s}\times \hcs s\times 
\cdots\times \hcs s.
$$
Moreover for any $s_0\in(\nu+3/2,s]$, there is $C>0$ such that for any
  $u_\ell\in\hcs s$, $\ell\in\{1,\ldots,k+1\}-\{j\}$, any $u_j\in\hcs{-s}$ 
\begin{equation}
\label{2.1.8}
\left|L(u_1,\ldots,u_{k+1})\right|\leq C \norma{L}_{\EL\nu
  N{k+1}}\norma{u_j}_{H^{-s}} \Big(\sum_{{1\leq\ell\leq k+1\atop \ell\not=j
  }}\norma{u_\ell}_{H^s}\prod_{{\ell'\not=\ell\atop \ell'\not=j}}
  \norma{u_{\ell'}}_{H^{s_0}}\Big). 
\end{equation}
\item[ii)] Let $k_1,k_2\in\Nn^*$, $\nu_1,\nu_2\in[0,+\infty)$, $1\leq 
  \ell\leq k_2+1$. For $M\in\EM{\nu_1}N{k_1}$, $L\in\EL{\nu_2}N{k_2+1}$ 
  with $N> 1+\max(\nu_1,\nu_2)$ define a $(k_1+k_2)$--linear form on 
  $\E^{k_1+k_2}$  
$$ 
(u_1,\ldots,u_{k_1+k_2})\to\tilde L(u_1,\ldots,u_{k_1+k_2}) 
$$  
by 
\begin{equation} 
\label{2.1.9} 
\tilde L(u_1,\ldots,u_{k_1+k_2})=
L(u_1,\ldots,u_{\ell-1},M(u_\ell,\ldots,u_{\ell+k_1-1}),
u_{\ell+k_1},\ldots,u_{k_1+k_2}). 
\end{equation} 
Then $\tilde L\in\EL{\nu_1+\nu_2+1}{N-\max(\nu_1,\nu_2)-1}{k_1+k_2}$ 
and the map $(M,L)\mapsto \tilde L$ is bounded from 
$\EM{\nu_1}N{k_1}\times \EL{\nu_2}N{k_2+1}$ to the preceding space. 
\end{itemize} 
\end{proposition} 
 
We shall denote, for any $N,\nu$ by 
\begin{equation} 
\label{2.1.14} 
\Sigma:\EL\nu N{k+1}\to \EM\nu N k 
\end{equation} 
the map given, using notation (\ref{2.1.6}), by 
$\Sigma(L)=M_{L,k+1}$. This is an isomorphism. 
 
In order to apply a Birkhoff procedure, it is necessary to verify that 
our framework is stable by Poisson brackets.
 
\begin{proposition}\label{prop2.1.4} 
    Let $k_1,k_2\in\Nn^*$, $\nu_1,\nu_2\in[0,+\infty)$, 
      $N>\frac{5}{2}+\max(\nu_1,\nu_2)$. Let $L_1\in\EL{\nu_1}N{k_1+1}$, 
      $L_2\in\EL{\nu_2}N{k_2+1}$, $\ell_1\in \{0,\ldots,k_1+1\}$, $\ell_2\in 
      \{0,\ldots,k_2+1\}$. Then 
      $\bigl\{\underline{L}\null_1^{\ell_1},\underline{L}\null_2^{\ell_2} 
      \bigr\}$ may be written  
    \begin{equation} 
    \label{2.1.11} 
    \bigl\{\underline{L}\null_1^{\ell_1},\underline{L}\null_2^{\ell_2}
    \bigr\}(u,\bar u) = \underline{L}\null_3^{\ell_1+\ell_2-1}(u,\bar
    u)
    \end{equation}   
    for a multilinear form $L_3\in 
    \EL{\nu_1+\nu_2+1}{N-\max(\nu_1,\nu_2)-1}{k_1+k_2}$. 
\end{proposition} 
 
\proof We can choose $s$ with 
$N-1>s>\frac{3}{2}+\max(\nu_1,\nu_2)$. By i) of proposition 
\ref{p.2.1.3}, $\underline{L}\null_i^{\ell_i}(u,\bar u)$ $i=1,2$ is then a 
smooth function on $\hcs s  $. Using \eqref{2.1.6} we may write for 
any $h\in\E$, $i=1,2$ 
\begin{equation*} 
\partial_u \underline{L}\null_i^{\ell_i}\cdot h  
= 
\sum_{j=1}^{\ell_i}L_i(u,\ldots,h,\ldots,u,\bar u,\ldots,\bar u) 
 =\sum_{j=1}^{\ell_i} \bigl\langle h, 
\underline{M}_{L_i,j}^{\ell_i-1}(u,\bar u)\bigr\rangle 
\end{equation*} 
where in the first sum $h$ stands at the $j$-th place. We have a 
similar formula for $\partial_{\bar u} \underline{L}\null_i^{\ell_i}. h 
$. In other words, we may write 
\begin{equation} 
\begin{split} 
\label{2.1.12} 
\nabla_u \underline{L}\null_i^{\ell_i}(u,\bar u)&= \sum_{j=1}^{\ell_i} 
\underline{M}_{L_i,j}^{\ell_i-1}(u,\bar u)  
\\ 
\nabla_{\bar u} \underline{L}\null_i^{\ell_i}(u,\bar u)&= 
\sum_{j=\ell_i+1}^{k_i+1}  
\underline{M}_{L_i,j}^{\ell_i}(u,\bar u).  
\end{split} 
\end{equation} 
By i) of proposition \ref{p.2.1.3.bis} these quantities are smooth 
functions of $u$ with values in $\hcs s  $, 
i.e. $\underline{L}\null_i^{\ell_i}\in   C_{s}^\infty(\hcs s,\C)$. We may thus 
apply definition \ref{d.1.2.2} and \eqref{1.2.18} to write 
\begin{equation} 
\begin{split} 
\label{2.1.12'} 
\bigl\{\underline{L}\null_1^{\ell_1},\underline{L}\null_2^{\ell_2}
\bigr\}(u,\bar u) = \im &\Big[
\sum_{j_2=1}^{\ell_2}\sum_{j_1=\ell_1+1}^{k_1+1}
L_2(u,\ldots,\underline{M}_{L_1,j_1}^{\ell_1}(u,\bar u),\ldots, u,\bar
u,\ldots,\bar u ) \\ &-
\sum_{j_2=\ell_2+1}^{k_2+1}\sum_{j_1=1}^{\ell_1} L_2(u,\ldots,u,\bar
u,\ldots,\underline{M}_{L_1,j_1}^{\ell_1-1}(u,\bar u),\ldots,\bar u )
\Big]
\end{split} 
\end{equation} 
where the $M$--term in the argument of $L_2$ stays at the $j_2$-th 
place. Since $M_{L_1,j_1}$ belongs to $\EM{\nu_1}{N}{k_1}$ we just 
have to apply (ii) of proposition \ref{p.2.1.3} to write this last 
expression in terms of a new multilinear form $L_3$. 
\qed 
 
\medskip 
 
In order to prove our main theorem we have to decompose the  
multilinear forms of $\EL\nu N{k+1}$ in the sum of a resonant  
and of a non-resonant part. 
 
\begin{definition} 
\label{d.2.1.4} 
(Non-resonant multilinear form) Fix $k\in \Nn$ and let $1\leq \ell\leq 
k+1$ be a fixed integer. 
\begin{itemize} 
\item If $2\ell\not =k+1$ we set $ \ELt\nu N{k+1}= \EL\nu N{k+1}$, 
$\EMt\nu N{k}= \EM\nu N{k}$ .
\item If $2\ell=k+1$ we define $\ELt\nu N{k+1}$ (resp. $\EMt\nu N{k}$) 
  as the subspace of those $L\in\EL\nu N{k+1} $ (resp. $M\in\EM\nu 
  N{k}$) such that respectively 
\begin{equation} 
\label{2.1.13} 
L(\Pi_{n_1}u_1,\ldots,\Pi_{n_{k+1}} u_{k+1})\equiv 0\ ,\quad 
\Pi_{n_{k+1}}M(\Pi_{n_1}u_1,\ldots,\Pi_{n_{k}} u_{k})\equiv 0 
\end{equation} 
for any $u_1,\ldots,u_{k+1}\in\E$ and any 
$(n_1,\ldots,n_{k+1})\in(\Nn^*)^{k+1}$ such that 
$$ 
\left\{n_1,\ldots,n_{\ell}\right\}=\left\{n_{\ell+1},\ldots,n_{k+1}\right\}. 
$$ 
\end{itemize}  
\end{definition}  
Remark that the map $\Sigma$ given by (\ref{2.1.14}) induces an isomorphism 
between $\ELt\nu N{k+1}$ and $\EMt \nu N{k}$. 
 
\begin{definition} 
\label{d.2.1.4.1} (Resonant multilinear form) 
Fix $k\in \Nn$ and let $1\leq\ell\leq k+1$. We define the space of
$\ell$--resonant multilinear forms
$\widehat{\mathcal{L}}_{k+1,\ell}^{\nu,N}$ as the subspace of those
$L\in\EL \nu N{k+1}$ verifying
\begin{equation} 
\label{2.1.13.1} 
L(\Pi_{n_1}u_1,\ldots,\Pi_{n_{k+1}} u_{k+1})\equiv 0\ , 
\end{equation} 
for any $u_1,\ldots,u_{k+1}\in\E $ and any 
$(n_1,\ldots,n_{k+1})\in(\Nn^*)^{k+1}$ such that 
$$ \left\{n_1,\ldots,n_{\ell}\right\}\not=\left\{n_{\ell+1},\ldots,
n_{k+1}\right\}.
$$ 
\end{definition} 
 
Remark that $\widehat{\mathcal{L}}_{k+1,\ell}^{\nu,N} = 0$ if $k$ is even
or $k$ is odd and $\ell \neq \frac{k+1}{2}$. If $k$ is odd and
$\ell = \frac{k+1}{2}$, one gets a direct sum decomposition
\begin{equation}
\label{decomp}
\EL \nu N{k+1} = \widehat{\mathcal{L}}_{k+1,\ell}^{\nu,N} \oplus
\tilde{\mathcal{L}}_{k+1,\ell}^{\nu,N}.
\end{equation}
The main feature of the above definitions is captured by the following 
proposition:
\begin{proposition} 
\label{p.jein} 
Assume that $L\in\EL\nu N{k+1}  $ is $\ell$ resonant. Then for any 
$a\in\Nn$, $a\geq 1$ one has 
$$ 
\big\{ \underline L^\ell,J_a \big\}\equiv 0.$$ 
\end{proposition} 
 
\proof Remark first that one has $2\ell=k+1$ and that $J_a(u,\bar 
u)=\langle \Pi_a u,\Pi_a\bar u\rangle$, from which, using 
\eqref{2.1.12'}, one gets 
$$ 
\big\{ \underline L^\ell,J_a \big\}=-\im\underline{\tilde L}^\ell
$$ 
with 
\begin{displaymath}
\begin{array}{l}
\displaystyle\tilde L(u_1,\ldots,u_{k+1})=\\  
\displaystyle\Big[\sum_{j=1}^{\ell}L(u_1,\ldots,{\Pi_a} u_j 
,\ldots,u_{k+1}) -\sum_{j=\ell+1}^{k+1}
L(u_1,\ldots,{\Pi_a} u_j,\ldots,u_{k+1}) \Big]. 
\end{array}
\end{displaymath}
 Then the above expression is equal to 
\begin{equation*} 
\begin{split}
\sum_{{n_1,\ldots,n_{k+1}}}\Big[ \sum_{j=1}^{\ell}
L(\Pi_{n_1}u_1,\ldots,\Pi_a\Pi_{n_j}u_j,\ldots,\Pi_{n_{k+1}}u_{k+1})
\makebox[2.5cm]{} \\ -\sum_{j=\ell+1}^{2\ell}
L(\Pi_{n_1}u_1,\ldots,\Pi_a\Pi_{n_j}u_j,\ldots,\Pi_{n_{k+1}}u_{k+1})\Big]
\\ =\sum_{n_1,\ldots,n_{k+1}} \Big[\sum_{j=1}^{\ell} \delta_{n_j,a}-
\sum_{j=\ell+1}^{2\ell}\delta_{n_j,a}\Big] L(\Pi_{n_1}u_1,\dots,
\Pi_{n_{2\ell}}u_{2\ell}).
\end{split}
\end{equation*} 
Since for an $\ell$ resonant form  
$$ 
\left\{n_1,\ldots,n_{\ell}\right\} =\left\{n_{\ell+1},
\ldots,n_{2\ell}\right\} ,
$$ 
the quantity $\sum_{j=1}^{\ell} 
\delta_{n_j,a}- \sum_{j=\ell+1}^{2\ell}\delta_{n_j,a}$ always vanishes.\qed 
 
\begin{definition} 
\label{2.1.d} 
For given integers $\ell, k$ satisfying $1\leq \ell\leq k+1$,  
we define an operator $\psi_\ell$ acting on $\EL\nu 
N{k+1}$ by 
\begin{equation} 
\begin{split} 
\label{2.1.10.a} 
\psi_\ell(L)(u_1,\ldots,u_{k+1}) \makebox[9cm]{}
\\  
=  
\Bigl[\sum_{j=1}^{\ell}L(u_1,\ldots,{\Lambda_m u}_j 
,\ldots,u_{k+1}) -\sum_{j=\ell+1}^{k+1}L(u_1,\ldots,{\Lambda_m 
  u}_j,\ldots,u_{k+1}) \Bigr]. 
\end{split} 
\end{equation} 
\end{definition} 
 
Remark that writing $G_{2}(u,\bar u) = \langle \Lambda_{m}u, \bar  
u\rangle$, and using \eqref{1.2.18} one gets 
\be \label{Gpsi} 
\big\{ \underline L^\ell,G_2\big\}(u,\bar u)= 
-\im \psi_\ell(L)(u,\ldots,u,\bar u,\ldots,\bar u) 
\ee 
where in the right hand side one has $\ell$ times $u$ and $k+1-\ell$  
times $\bar u$.  
 
\begin{proposition} 
\label{prop3.5} 
There is a zero measure subset $\N$ of $(0,+\infty)$ such that for any 
$k\in\Nn^*$, any $m\in (0,+\infty)-\N$, any $0\leq \ell\leq k+1$, there is 
a $\bar \nu\in\R_+$, and for any $(\nu,N)\in\R_+\times \Nn, N>2$, there is 
an operator 
\begin{equation} 
\label{2.1.18} 
\psi_\ell^{-1}:\ELt\nu N{k+1}\to \ELt{\nu+\bar\nu} N{k+1} 
\end{equation}    
such that for any $L\in \ELt\nu N{k+1}$, $\psi_\ell(\psi_\ell^{-1}(L))=L$.  
Moreover there exists $C>0$ such that 
\begin{equation} 
\label{2.1.18.1} 
\norma{\psi_\ell^{-1}(L)}_{\EL{\nu+\bar\nu}N{k+1}} \leq C 
\norma{L}_{\EL{\nu}N{k+1}}.  
\end{equation}  
\end{proposition} 
\proof We reduce the proof to  proposition 2.2.4 of
\cite{DS2}. Let
$\rho:\left\{1,\ldots,k+1\right\}\to\left\{-1,1\right\}$ be the map given
by $\rho(j)=1$ if $j=1,\ldots,\ell$ and $\rho(j)=-1$ if
$j=\ell+1,\ldots,k+1$, and for $M\in\EM \nu N k$ define
\begin{equation} 
\label{2.1.18.b} 
\begin{array}{l}
\displaystyle\tilde\psi_\ell(M)(u_1,\ldots,u_{k})=\\
\displaystyle\sum_{j=1}^{k}\rho(j)
M(u_1,\ldots,\Lambda_mu_j,\ldots,u_{k})+\rho(k+1)\Lambda_mM
(u_1,\ldots,u_{k}).
\end{array}
\end{equation} 
One has, if $\Sigma$ is the map defined in \eqref{2.1.14}, 
\begin{equation} 
\label{2.1.18.a} 
\Sigma^{-1}\circ\tilde\psi_\ell(M)=\psi_\ell\circ\Sigma^{-1}(M)
\end{equation} 
for any $M \in \mathcal{M}^{\nu,N}_k$ such that $\tilde\psi_\ell(M)$
belongs to $\mathcal{M}^{\nu',N}_k$ for some $\nu'\geq 0$.
By proposition 2.2.4 in \cite{DS2}, there are $\bar\nu\in\R_+$ and an
operator $\tilde \psi_\ell^{-1}:\EMt\nu Nk\to \EMt{\nu+\bar\nu} Nk$
such that for any $M\in \EMt\nu Nk$,
$\tilde\psi_\ell(\tilde\psi_\ell^{-1}(M))=M$ and such that the
equivalent for $M$ of the estimate (\ref{2.1.18.1}) holds true. We
just set $\psi_\ell^{-1}=\Sigma^{-1}\circ \tilde \psi_\ell^{-1}
\circ\Sigma$, and the conclusion follows from equation
(\ref{2.1.18.a}).  \qed
\medskip
 
The construction of the operator $\tilde \psi_\ell^{-1}$ in \cite{DS2} 
relies in an essential way on the spectral assumption (\ref{Kn}) and 
(\ref{115}), i.e. on the fact that $M$ is a Zoll manifold. For the  
reader's convenience, we give a direct proof of  
proposition \ref{prop3.5} in the case where $M={\mathbb{S}}^d$ and $V=0$.  
In this case, the eigenvalues $\lambda_n$ of $P$ and 
$\omega_n$ of $\Lambda_m$ are respectively given by 
\begin{equation} 
\label{A.3} 
\lambda_n=\sqrt{n(n+d-1)}\ ,\quad \omega_n=\sqrt{\lambda_n^2+m^2}\ , 
\end{equation} 
and moreover $P\Pi_n=\lambda_n\Pi_n$, 
$\Lambda_m\Pi_n=\omega_n\Pi_n$. Thus, from equation (\ref{2.1.10.a}) one 
has 
\begin{equation} 
\begin{split} 
\label{A.1} 
&\psi_\ell(L)(u_1,\ldots,u_{k+1}) \\ &= \sum_{n_1,\ldots,n_{k+1}}
 (\omega_{n_1}+\cdots+\omega_{n_\ell}-\omega_{n_{\ell+1}}-\cdots
 -\omega_{n_{k+1}} ) L(\Pi_{n_1}u_1,\ldots,\Pi_{n_{k+1}}u_{k+1}).
\end{split} 
\end{equation} 
Remark also that, if $L\in\ELt \nu N{k+1}$, then the sum is restricted 
to those $(n_1,\ldots,n_{k+1})$ such that
$$ \left\{n_1,\ldots,n_{\ell}\right\}\not=
\left\{n_{\ell+1},\ldots,n_{k+1}\right\}\ .
$$ 
The following proposition was proved in \cite{DS1} (see Proposition 
4.8) and is also a minor variant of theorem 3.12 of \cite{BG}. 
 
\begin{proposition} 
\label{A.4} 
There is a zero measure subset $\N$ of $(0,+\infty)$ such that for any 
$m\in (0,+\infty)-\N$ and any $k\in\Nn^*$, there are $ c>0$ and 
$\bar \nu\in\R_+$ such that for any $0\leq \ell\leq k+1$, one has 
\begin{equation} 
\label{A.2} 
\left|\omega_{n_1}+\cdots+\omega_{n_\ell}-\omega_{n_{\ell+1}}-\cdots 
-\omega_{n_{k+1}} \right|\geq c\mu (n_1,\ldots,n_{k+1})^{-\bar\nu}  
\end{equation} 
for any choice of $(n_1,\ldots,n_{k+1})$ such that 
$$ \left\{n_1,\ldots,n_{\ell}\right\}\not=
\left\{n_{\ell+1},\ldots,n_{k+1}\right\}\ .
$$ 
\end{proposition}  
It is now immediate to obtain the 
 
\noindent 
{\bf Proof of Proposition \ref{prop3.5} in the case $M={\mathbb{S}}^d$, 
  $V\equiv 0$.} Given $L\in\ELt \nu N{k+1}$ define 
\begin{equation}
\label{A.12} 
\tilde L(u_1,\ldots,u_{k+1}) = \sum_{n_1,\ldots,n_{k+1}}
 \frac{L(\Pi_{n_1}u_1,\ldots,\Pi_{n_{k+1}}u_{k+1})}{
 (\omega_{n_1}+\cdots+\omega_{n_\ell}-\omega_{n_{\ell+1}}-\cdots
 -\omega_{n_{k+1}} ) }.
\end{equation}  
Then by (\ref{A.2}) one has $\tilde L\in\ELt {\nu+\bar\nu} N{k+1}$, 
and by (\ref{A.12}) $\psi_\ell(\tilde L)=L$; finally also the estimate 
(\ref{2.1.18.1}) immediately follows. On a general Zoll manifold, the 
construction of the map $L \to \tilde{L}$ is made in~\cite{DS1} through an
approximation argument and a suitable use of Neumann series. \qed 

\medskip

Finally we end this subsection with two lemmas that will be useful to 
verify that certain Hamiltonian functions are real valued.
\begin{lemma} 
\label{lemma3.13} 
Assume $m\in (0,+\infty)-\N$ and let $L\in \ELt\nu N{k+1}$.\\ i)
Assume that for any $u\in\E$, $\psi_\ell(L)(u,\ldots,u,\bar
u,\ldots,\bar u)=0$ (where one has $\ell$ times $u$ and $k+1-\ell$
times $\bar u$). Then $\underline L^\ell(u,\bar u)=0$.\\ ii) Assume
$\bigl\{\Im \underline{L}^\ell,G_2\bigr\} \equiv 0$. Then
$\Im\underline{L}^\ell(u,\bar{u}) \equiv 0$.
\end{lemma} 

\proof i) Let $\mathfrak{S}_{\ell,k}$ be the product of the group of
permutations of $\{1,\ldots,\ell\}$ by the group of permutations of
$\{\ell+1,\ldots,k+1\}$. For
$(\sigma,\sigma')\in\mathfrak{S}_{\ell,k}$ define
\[((\sigma,\sigma')\cdot L)(u_1,\ldots,u_{k+1})=L(u_{\sigma(1)},
\ldots,u_{\sigma(\ell)},u_{\sigma'(\ell+1)},
\ldots,u_{\sigma'(k+1)}).\] 
Replacing $L$ by
\[\frac{1}{\ell!(k+1-\ell)!} \sum_{a\in\mathfrak{S}_{\ell,k}}(a\cdot L)
(u_1,\ldots,u_{k+1})\] 
does no affect the hypotheses nor the conclusion (since $\psi_\ell$
commutes to the $\mathfrak{S}_{\ell,k}$-action), so we can assume that
$L$ -- and thus $\psi_\ell(L)$ -- is
$\mathfrak{S}_{\ell,k}$-invariant. Write the assumption
$\psi_\ell(L)(u,\ldots,u,\bar u,\ldots,\bar u)=0$ with
\[u=u_1+\cdots+u_\ell+\overline{u_{\ell+1}}+\cdots+\overline{u_{k+1}}\]
for arbitrary $u_j$'s belonging  
to $\mathcal{E}$. If one expands this expression by multilinearity,
sorts the different contributions according to their homogeneity
degree in $u_j,\bar{u}_j$, and uses the
$\mathfrak{S}_{\ell,k}$-invariance, one gets
\begin{equation}\label{inj1}
\psi_\ell(L)(u_1,\ldots,u_{k+1})=0
\end{equation}
for any $u_1,\ldots,u_{k+1}$ in $\mathcal{E}$. Take a family of
positive integers $(n_1,\ldots,n_{k+1})$ such that
$\{n_1,\ldots,n_\ell\}\neq\{n_{\ell+1},\ldots,n_{k+1}\}$ if
$2\ell=k+1$. We apply \eqref{inj1} taking for all $u_j$ an
eigenfunction associated to an eigenvalue $\lambda_{n_j}\in K_{n_j}$
$j=1,\ldots,k+1$ so that $\Lambda_mu_j=\omega_{n_j}u_j$,
$\omega_{n_j}= \sqrt{m^2+\lambda_{n_j}^2}$.
By \eqref{2.1.10.a} we obtain
\[\Big(\sum_{j=1}^\ell
\omega_{n_j}-
\sum_{\ell+1}^{k+1}\omega_{n_j}\Big)L(u_1,\ldots,
u_{k+1})=0.\] By proposition 2.2.1 and formula (2.2.3) of \cite{DS2}
(see also proposition \ref{A.4} of the present paper in the case of
the sphere), the first factor is nonzero for $m\in (0,+\infty)-\N$, so
$L(u_1,\ldots,u_{k+1})=0$ for any family $(u_1,\ldots,u_{k+1})$ of the
preceding form. The definition of $\ELt\nu N{k+1}$ implies that
$\underline L^\ell(u,\bar u)=0$.

ii) We may write when $\ell \neq \frac{k+1}{2}$
$\Im\underline{L}^\ell(u,\bar{u}) =
\underline{\Gamma}_1^\ell(u,\bar{u}) +
\underline{\Gamma}_2^{k+1-\ell}(u,\bar{u})$ for $\Gamma_1 \in \ELt\nu
N{k+1}$, $\Gamma_2 \in \tilde{\mathcal{L}}^{\nu,N}_{k+1,k+1-\ell}$. By
homogeneity, $\bigl\{G_2,\underline{\Gamma}_1^\ell +
\underline{\Gamma}_2^{k+1-\ell}\bigr\} \equiv 0$  implies that
$\bigl\{G_2,\underline{\Gamma}_1^\ell\bigr\} =$
\mbox{$\bigl\{G_2,\underline{\Gamma}_2^{k+1-\ell}\bigr\} \equiv 0$},
whence $\underline{\Gamma}_1^\ell = \underline{\Gamma}_2^{k+1-\ell}
=0$ by \eqref{Gpsi} and assertion i). If $\ell = \frac{k+1}{2}$, we
have $\Im\underline{L}^\ell(u,\bar{u}) =
\underline{\Gamma}^\ell(u,\bar{u})$ for a $\Gamma \in
\widetilde{\mathcal{L}}^{\nu,N}_{k+1}$, and the result follows again
from i).  \qed

\begin{lemma} 
\label{lemma3.13bis} 
Assume $m\in (0,+\infty)-\N$ and $k$ odd. Set $\ell = \frac{k+1}{2}$
and consider $L_1\in \ELt\nu N{k+1}$ and $L_2 \in
\widehat{\mathcal{L}}_{k+1,\ell}^{\nu,N}$.  Set $L = L_1+L_2$ and
assume that for any $u\in\E$, $\underline L^\ell(u,\bar u)$ is real
valued. Then $\underline {L_1}^\ell(u,\bar u)$ and $\underline
{L_2}^\ell(u,\bar u)$ are real valued.
\end{lemma}  

\proof Since $\underline L^\ell(u,\bar u)$ is real valued,
$\{\Im\underline L^\ell,G_2\}(u,\bar u)=0$. As  \[\{\underline
L^\ell,G_2\}(u,\bar u) =\{\underline {L_1}^\ell,G_2\}(u,\bar u)\] by
proposition \ref{p.jein}, this yields $\{\Im\underline
{L_1}^\ell,G_2\}(u,\bar u)=0$. Now, ii) of lemma \ref{lemma3.13}
implies $\Im\underline {L_1}^\ell(u,\bar u)=0$. Therefore,
$\underline{L_1}^\ell(u,\bar u)$ and $\underline {L_2}^\ell(u,\bar u)$
are real valued.  \qed

\subsection{Proof of theorem \ref{thm2}.}\label{birk} 

We use a Birkhoff scheme to put the Hamiltonian system with the
Hamiltonian $G$ of \eqref{1.2.12} in normal form. Having fixed some
$r_0\geq 1$, the idea is to construct iteratively for $r =
0,\ldots,r_0$, $\U_r$ a neighborhood of $0$ in $H^s(M,\C)$ for $s\gg 
1$, a canonical transformation $\Tr_r$, defined on $\U_r$, an
increasing sequence $(\nu_r)_{r=1,\ldots,r_0}$ of positive numbers,
and functions $\Ze^{(r)}, P^{(r)}, \resto^{(r)}$ such that
\begin{equation} 
\label{b.1} 
G^{(r)}:=G\circ\Tr_r=G_2+\Ze^{(r)}+P^{(r)}+\resto^{(r)}. 
\end{equation} 
Moreover, these functions will decompose as 
\begin{eqnarray} 
\label{b.2} 
\Ze^{(r)}&=&\sum_{j=1}^r \Ye_j 
\\ 
\label{b.3} 
P^{(r)}&=&\sum_{j=r+1}^{r_0} Q^{(r)}_j 
\end{eqnarray} 
where $\Ye_j$ is in $\mathcal{H}_{j+2}^s(\nu_j)$ and Poisson commutes
with $J_n$ for any $n$, $Q^{(r)}_j$ is in
$\mathcal{H}_{j+2}^s(\nu_r)$, by convention $P^{(r_{0})}=0$, and
$\resto^{(r)}\in C_s^\infty(\U_r,\R)$ has a zero of order $r_0+3$ at the
origin.
 
First remark that the Hamiltonian (\ref{1.2.12}) has the form 
(\ref{b.1}), (\ref{b.2}), (\ref{b.3}) with $r=0$ and $\Tr_r=I$,  
$P^{(0)}$ being the Taylor's polynomial of $\tilde{G}$ at degree 
$r_{0}$ (see lemma \ref{lem:G}). We 
show now how to pass from $r$ to $r+1$ provided one is able to solve
the homological equation below. 
 
\begin{lemma} 
\label{lemma3.14} 
Assume we are given $0<\nu_r$ and functions $\Ze^{(r)}, P^{(r)}, 
 \resto^{(r)}$ satisfying the above conditions. Assume that there are
 $\nu'_r>\nu_r$ and  a function $F^{(r+1)}$ of $(u,\bar u)$
with the properties that   
\begin{eqnarray} 
\label{b.4} 
F^{(r+1)}&\in& \mathcal{H}_{r+3}^s(\nu'_r)
\\ 
\label{b.5} 
\{F^{(r+1)},G_2\}&\in& \mathcal{H}_{r+3}^s(\nu'_r). 
\end{eqnarray} 
Assume moreover one is able to choose $F^{(r+1)}$ with the further
property that $\Ye_{r+1}$ defined by
  \begin{equation} 
\label{b.10} 
\Ye_{r+1}=\bigl\{F^{(r+1)},G_2\bigr\}+Q^{(r)}_{r+1} 
\end{equation} 
Poisson commutes with $J_n$ for any $n$.  Denote by $\Phi^t_{r+1}$ the
flow generated by $X_{F^{(r+1)}}$. Then, there are $\nu_{r+1}>\nu'_r$
and, for large enough $s$, a sufficiently small neighborhood
$\U_{r+1}$ of the origin of $H^s(M,\C)$, such that $G^{(r+1)}=
G^{(r)}\circ\Phi^1_{r+1}$ has the same structure as $G^{(r)}$ but with
$r$ replaced by $r+1$ and $\U_r$ replaced by $\U_{r+1}$.
\end{lemma}  

\proof If 
$\U_{r+1}$ is a sufficiently small neighborhood of the origin of 
$H^s(M,\C)$, then $\Phi_{r+1}^1:\U_{r+1}\to\U_{r}$ is well defined. 
We decompose $G^{(r)}\circ \Phi^1_{r+1}$ as follows 
\begin{eqnarray} 
\label{b.6} 
G^{(r)}\circ \Phi^1_{r+1}&=& G_2+\bigl\{F^{(r+1)},G_2\bigr\}+Q^{(r)}_{r+1}+ 
\Ze^{(r)} + \resto^{(r)}\circ\Phi^1_{r+1} 
\\ 
\label{b.7} 
&+& P^{(r)}\circ\Phi^1_{r+1}-Q^{(r)}_{r+1} 
\\ 
\label{b.8} 
&+& \Ze^{(r)}\circ\Phi^1_{r+1}-\Ze^{(r)} 
\\ 
\label{b.9} 
&+& G_2\circ\Phi^1_{r+1}-G_2-\bigl\{F^{(r+1)},G_2\bigr\}.
\end{eqnarray} 
Using the fact that $\Ye_{r+1}$ Poisson commutes with $J_n$ for any
$n$ and belongs to $\mathcal{H}_{r+3}^s(\nu'_r) \subset
\mathcal{H}_{r+3}^s(\nu_{r+1})$ by \eqref{b.5}, we may define
$\Ze^{(r+1)}:=\Ze^{(r)}+\Ye_{r+1}$.

If $s$ is large enough, \eqref{b.4} implies that $F^{(r+1)} \in
C^\infty_s(\U_r,\R)$, and we may apply lemma \ref{l.1.2.3} with
$F=F^{(r+1)}$ to $P^{(r)}\circ\Phi^1_{r+1}$ and
$\Ze^{(r)}\circ\Phi^1_{r+1}$.  Using proposition~\ref{prop2.1.4} to
write the iterated Poisson brackets of the right hand side
of~\eqref{1.2.23} in terms of multilinear forms, we thus see that
\eqref{b.7}, \eqref{b.8} may be decomposed in a sum of elements of
$\mathcal{H}^s_{j+2}(\nu_{r+1})$ for $s, \nu_{r+1}$ large enough and
$j= r+2,\ldots,r_0$. Consequently these two terms will contribute to
$P^{(r+1)}, \resto^{(r+1)}$ in \eqref{b.1} written with $r$ replaced
by $r+1$. In the same way, lemma~\ref{l.1.2.3} applied to
$G_2\circ\Phi_{r+1}^1$ shows that, for large enough $r$ and
$\nu_{r+1}$, \eqref{b.9} gives a contribution to $P^{(r+1)}+
\resto^{(r+1)}$ in \eqref{b.1} at step $r+1$. The conclusion follows.
\qed
\medskip

Let us remark that the above lemma implies theorem 2.6. Actually, if
we are able to apply lemma 3.17 up to step $r_0 -1$, we get (3.31)
with $r=r_0$, which is the conclusion of the theorem. Our remaining
task is thus to solve the homological equation (3.36). This will be
achieved in the following lemma.

\begin{lemma}
\label{hom.true}
Let $2\leq r\leq r_0+2, \nu\in \R_+^*$, and assume $m\in
(0,+\infty)-\N$. For any $Q\in \Hs{r+1}(\nu)$ there are $\nu'>\nu$,
$F\in\Hs{r+1}(\nu')$ and $\Ye\in\Hs{r+1} (\nu)$, with $\Ye$ which
Poisson commutes with $J_n$ for any $n\geq1$, such that
\begin{equation}
\label{hom.1}
\left\{F,G_2\right\}+Q=\Ye\, .
\end{equation}
As a consequence one also has $\{F,G_2\}\in\Hs{r+1}(\nu)$.
\end{lemma}
\proof If $r+1$ is odd then we define $\Ye=0$. As $Q$ is
in $\mathcal{H}_{r+1}^s(\nu)$, it decomposes in the form
\begin{equation} 
\label{F.11} 
Q=\sum_{\ell=0}^{r+1} \underline{L_{\ell}}\null^\ell 
\end{equation} 
where $L_{\ell}$ are multilinear forms in $\EL{\nu} {+\infty}{r+1}$.  
We remark that, since $r+1$ is odd, the $L_{\ell}$ are all  
non-resonant, i.e. $L_{\ell} \in \ELt{\nu} {+\infty}{r+1}$. Therefore  
by proposition \ref{prop3.5}, we can define  
$F_{\ell}\in \ELt{\nu+\bar\nu} {+\infty}{r+1}$ by 
\be\label{F1} 
F_{\ell}=-\im \Psi_{\ell}^{-1}(L_{\ell}) 
\ee  
and in view of  \eqref{Gpsi}, the Hamiltonian function 
\be\label{F2} 
F=\sum_{\ell=0}^{r+1} \underline{F_{\ell}}\null^\ell 
\ee 
satisfies the homological equation \eqref{hom.1}.  

If $r+1$ is even, set $\tilde{L}_{\ell} = L_{\ell}
$ if $\ell \neq \frac{r+1}{2}$. When $\ell=\frac{r+1}{2}$, write  
\[ L_{\frac{r+1}{2}} = Y + 
\tilde{L}_{\frac{r+1}{2}} \in 
\widehat{\mathcal{L}}^{\nu,+\infty}_{r+1,\frac{r+1}{2}} \oplus  
\widetilde{\mathcal{L}}^{\nu,+\infty}_{r+1,\frac{r+1}{2}}\] 
using decomposition \eqref{decomp}. Then if
$\Ye:=\underline{Y}^{\frac{r+1}{2}}$,  
\[Q-\Ye=\sum_{\ell=0}^{r+1} 
\underline{(\tilde{L}_{\ell})}\null^\ell \] 
and if we define  $F$ by    
  \eqref{F2} with $F_{\ell} = 
  -\im\psi_\ell^{-1}(\tilde{L}_{\ell})$, we still obtain
  that equation \eqref{hom.1} is satisfied. 

It remains to show that $F$ is real valued. 
As $Q$ is real, using \eqref{F.11} yields for any
$\ell\in\{0,\ldots,r+1\}$ 
\begin{equation} 
\label{F.13} 
\overline{\underline{L_{\ell}}\null^\ell}
=\underline{L_{r+1-\ell}}\null^{r+1-\ell} 
\end{equation} 
by homogeneity. If $r+1$ is even, \eqref{F.13} implies that  
$\underline{L_{\frac{r+1}{2}}}\null^\frac{r+1}{2}$ is real
valued. Using lemma \ref{lemma3.13bis}, we obtain that $\Ye$ is real
valued (remark that if $r+1$ is odd, $\Ye=0$ is also real  
valued). Therefore, $\left\{F,G_2\right\}$ is real valued by  
\eqref{hom.1}. So $\left\{\Im F,G_2\right\}=0$ which implies by  
homogeneity that $\bigl\{\Im 
\underline{F_{\ell}}\null^\ell,G_2\bigr\}=0$ for any $\ell$.  We may now use 
lemma 
\ref{lemma3.13} to obtain that $\Im 
\underline{F_{\ell}}\null^\ell=0$ for any
$\ell\in\{0,\ldots,r+3\}$.  Therefore, $F$ is real valued.\qed

\subsection{Proof of theorem \ref{thm1}.}\label{sec:proofthm1}   
 
Let $\Tr$ be the canonical transformation defined in theorem \ref{thm2}. 
Define on $\U=\Tr(\V)$ the function 
$$
E(u,\bar u):=\sum_{n\geq 1}n^{2s}J_n\circ\Tr^{-1}(u,\bar u).
$$ We shall control $E(u,\bar{u})$ along long time intervals. To take
into account the loss of derivatives coming from the linear part of
the equation, we proceed by regularization. Fix $\sigma=s+1$ and take
the Cauchy data such that  $u_0=\epsilon (\Lambda_m^{-1/2}v_{1} + \im
\Lambda_m^{1/2}v_{0})/\sqrt{2}$ is in $H^\sigma(M,\C)\cap\U$. Let
$u(t)\equiv u(t,.)$ be the corresponding solution of $\dot
u=X_G(u)\equiv \im \nabla_{\bar{u}} G(u,\bar{u})$. Since $X_G$ is semilinear and
$H^\sigma$ is its domain, as far as $\norma{u(t)}_{H^s}<\infty$ one
has $u(t)\in H^\sigma$. Thus, as far as $u(t)\in\U$
\begin{equation}\label{F.17}
\frac{dE}{dt}=\partial E\cdot X_G=\left\{G,E\right\}
\end{equation}
which is well defined since $E\in C^\infty(\U,\R)$, with $\U\subset
H^s$ and $X_G(u)\in H^s$ for $u\in H^\sigma$. So we may write 
\[\left\{G,E\right\}(u,\bar u)  =\displaystyle\sum_{n\geq
1}n^{2s}\left\{G,J_n\circ\Tr^{-1}\right\}.\]
If we use \eqref{1.2.21}, \eqref{for.nor} and \eqref{stime1} we get then
\begin{equation}\label{F.18}
\begin{array}{ll}
\left\{G,E\right\} & 
=\displaystyle\sum_{n\geq
1}n^{2s}\left\{G\circ\Tr,J_n\right\}\circ\Tr^{-1}\\ &
=\displaystyle\sum_{n\geq
1}n^{2s}\left\{G_2+\Ze+\resto,J_n\right\}\circ\Tr^{-1}\\ &
=\left\{\resto\circ\Tr^{-1},E\right\}.
\end{array}
\end{equation}
Thus
\[\frac{dE}{dt}=\left\{\resto\circ\Tr^{-1},E\right\}\]
which, by taking an approximating sequence is seen to hold also for
initial data which are not in $H^\sigma$, but only in $\U$.
Using \eqref{stime2} one has
\begin{equation}\label{F.14}
\left|\frac{dE(t)}{dt}\right|\leq C\Vert{u(t,\cdot )}\Vert_{H^s}^{r+3}.
\end{equation}
Remark that by definition of $E(u,\bar{u})$ and because $\Tr(0) = 0$,
as long as $u$ stays in a small enough neighborhood of 0, we have
\begin{equation}
\label{la.e}
\frac{1}{2}E(u,\bar u)\leq \norma{u}^2_{H^s}\leq 2E(u,\bar u).
\end{equation}
We deduce then by integration of \eqref{F.14} the estimate
\begin{equation}\label{F.16}
\Vert{u(t,\cdot )}\Vert_{H^s}^2 \leq C'\left(\Vert{u_0}\Vert_{H^s}^2
+\left|\int_0^t\Vert{u(\tau,\cdot )}\Vert_{H^s}^{r+3}d\tau\right|\right)
\end{equation}
which holds true as long as $u$ remains in a small enough neighborhood
of 0. It is classical to deduce from this inequality that there are
$C>0, c>0, \epsilon_0>0$ such that, if the
Cauchy data $u_0$ is in the $H^s$ ball of center 0 and radius $\epsilon
<\epsilon_0$, the solution exists over an interval of length at least
$c\epsilon^{-r-1}$, and for any $t$ in that interval $\Vert
u(t,\cdot )\Vert_{H^s} \leq C\epsilon$. This concludes the proof.
\qed

\begin{remark} 
The proof of \eqref{estimJn} is similar. As in \eqref{F.17} and
\eqref{F.18}, we see that
\[\frac{dJ_n\circ \Tr^{-1}(u,\bar u)}{dt}=\left\{\resto\circ\Tr^{-1},
J_n\right\}(u,\bar u)\] 
which together with the bound $\Vert{u(t,\cdot 
)}\Vert_{H^s}\leq C_1\epsilon$ yields
\begin{equation}\label{F.19}
|J_n\circ\Tr^{-1}(u(t),\bar{u}(t))-J_n\circ\Tr^{-1}(u_0,\bar{u}_0)|\leq
 \frac{C\epsilon^3}{n^{2s}}
\end{equation}
for times $|t|\leq\epsilon^{-r}$. Finally, using \eqref{F.19},
\eqref{stime} and the inequality
\begin{displaymath}
\begin{array}{l}
|J_n(u(t),\bar{u}(t))-J_n(u_0,\bar{u}_0)|\leq
 |J_n(u(t),\bar{u}(t))-J_n\circ\Tr^{-1}(u(t),\bar{u}(t))|\\
\hspace{0.4cm}+|J_n\circ\Tr^{-1}(u(t),\bar{u}(t))-J_n \circ\Tr^{-1}(u_0,
\bar{u}_0)|+|J_n(u_0,\bar{u}_0)-J_n\circ\Tr^{-1}(u_0,\bar{u}_0)|
\end{array}
\end{displaymath}
implies \eqref{estimJn}. 
\end{remark}


\begin{thebibliography}{10}

\bibitem{Bam93}
D.~Bambusi and A.~Giorgilli.
\newblock Exponential stability of states close to resonance in
  infinite-dimensional {H}amiltonian systems.
\newblock {\em J. Statist. Phys.}, 71(3-4):569--606, 1993.

\bibitem{BG}
D.~Bambusi and B.~Gr{\'e}bert.
\newblock Birkhoff normal form for pdes with tame modulus.
\newblock {\em Duke Math. J.}, To appear.

\bibitem{BN98}
D.~Bambusi and N.~N. Nekhoroshev.
\newblock A property of exponential stability in the nonlinear wave equation
  close to main linear mode.
\newblock {\em Physica D}, 122:73--104, 1998.

\bibitem{Bam03}
D.~Bambusi.
\newblock Birkhoff normal form for some nonlinear {PDE}s.
\newblock {\em Comm. Math. Physics}, 234:253--283, 2003.

\bibitem{Bo96}
J.~Bourgain.
\newblock Construction of approximative and almost-periodic solutions of
  perturbed linear {S}chr{\"o}dinger and wave equations.
\newblock {\em Geometric and Functional Analysis}, 6:201--230, 1996.

\bibitem{Bo96b}
J.~Bourgain.
\newblock On the growth in time of higher Sobolev norms of smooth
solutions of Hamiltonian PDE.  
\newblock {\em Internat. Math. Res. Notices}, no. 6, 277--304, 1996.

\bibitem{Bo04}
J.~Bourgain.
\newblock Remarks on stability and diffusion in high-dimensional
              {H}amiltonian systems and partial differential equations.
  \newblock { \em Ergodic Theory Dynam. Systems}, 24:1331--1357, 2004.

\bibitem{CheMa}
P.~R. Chernoff and J.~E. Marsden.
\newblock {\em Properties of infinite dimensional {H}amiltonian systems}.
\newblock Springer-Verlag, Berlin, 1974.
\newblock Lecture Notes in Mathematics, Vol. 425.

\bibitem{CV}
Y.~Colin~de Verdi{\`e}re.
\newblock Sur le spectre des op{\'e}rateurs elliptiques {\`a} bicaract{\'e}ristiques
  toutes p{\'e}riodiques.
\newblock {\em Comment. Math. Helv.}, 54(3):508--522, 1979.

\bibitem{CW}
W.~Craig and C.~E. Wayne.
\newblock Newton's method and periodic solutions of nonlinear wave equations.
\newblock {\em Comm. Pure Appl. Math.}, 46:1409--1498, 1993.

\bibitem{D1}
J.-M. Delort.
\newblock Temps d'existence pour l'{\'e}quation de {K}lein-{G}ordon
  semi-lin{\'e}aire {\`a}\ donn{\'e}es petites p{\'e}riodiques.
\newblock {\em Amer. J. Math.}, 120(3):663--689, 1998.

\bibitem{D}
J.-M. Delort.
\newblock Existence globale et comportement asymptotique pour l'{\'e}quation de
  {K}lein-{G}ordon quasi lin{\'e}aire {\`a} donn{\'e}es petites en dimension 1.
\newblock {\em Ann. Sci. {\'E}cole Norm. Sup. (4)}, 34(1):1--61, 2001.

\bibitem{DS1}
J.-M. Delort and J.~Szeftel.
\newblock Long--time existence for small data nonlinear {K}lein--{G}ordon
  equations on tori and spheres.
\newblock {\em Internat. Math. Res. Notices}, 37:1897--1966, 2004.

\bibitem{DS2}
J.-M. Delort and J.~Szeftel.
\newblock Long--time existence for semi--linear {K}lein--{G}ordon equations
  with small Cauchy data on {Z}oll manifolds.
\newblock {\em Preprint}, 2004.

\bibitem{DS3}
J.-M. Delort and J.~Szeftel.
\newblock Bounded almost global solutions for non Hamiltonian semi-linear
  {K}lein--{G}ordon equations with radial data on compact revolution
  hypersurfaces.
\newblock {\em Preprint}, 2005.

\bibitem{DG}
J.~J. Duistermaat and V.~W. Guillemin.
\newblock The spectrum of positive elliptic operators and periodic
  bicharacteristics.
\newblock {\em Invent. Math.}, 29(1):39--79, 1975.

\bibitem{G}
V.~Guillemin.
\newblock Lectures on spectral theory of elliptic operators.
\newblock {\em Duke Math. J.}, 44(3):485--517, 1977.

\bibitem{K}
S.~Klainerman.
\newblock The null condition and global existence to nonlinear wave equations.
\newblock In {\em Nonlinear systems of partial differential equations in
  applied mathematics, Part 1 (Santa Fe, N.M., 1984)}, volume~23 of {\em
  Lectures in Appl. Math.}, pages 293--326. Amer. Math. Soc., Providence, RI,
  1986.

\bibitem{K1}
S.~B. Kuksin.
\newblock {\em Nearly integrable infinite-dimensional {H}amiltonian Systems}.
\newblock Springer-Verlag, Berlin, 1993.

  
\bibitem{KP96}
S.~B. Kuksin and J.~P{\"o}schel.
\newblock Invariant {C}antor manifolds of quasi-periodic oscillations for a
  nonlinear {S}chr{\"o}dinger equation.
\newblock {\em Ann. Math.}, 143:149--179, 1996.


\bibitem{MTT}
K.~Moriyama, S.~Tonegawa, and Y.~Tsutsumi.
\newblock Almost global existence of solutions for the quadratic semilinear
  {K}lein-{G}ordon equation in one space dimension.
\newblock {\em Funkcial. Ekvac.}, 40(2):313--333, 1997.

\bibitem{Muj}
J.~Mujica.
\newblock {\em Complex analysis in {B}anach spaces}, volume 120 of {\em
  North-Holland Mathematics Studies}.
\newblock North-Holland Publishing Co., Amsterdam, 1986.
\newblock Holomorphic functions and domains of holomorphy in finite and
  infinite dimensions, Notas de Matem{\'a}tica [Mathematical Notes], 107.

\bibitem{Nik86}
N.~Nikolenko.
\newblock The method of {P}oincar{\'e} normal form in problems of integrability
  of equations of evolution type.
\newblock {\em Russ. Math. Surveys}, 41:63--114, 1986.

\bibitem{OTT}
T.~Ozawa, K.~Tsutaya, and Y.~Tsutsumi.
\newblock Global existence and asymptotic behavior of solutions for the
  {K}lein-{G}ordon equations with quadratic nonlinearity in two space
  dimensions.
\newblock {\em Math. Z.}, 222(3):341--362, 1996.

\bibitem{Sh}
J.~Shatah.
\newblock Normal forms and quadratic nonlinear {K}lein-{G}ordon equations.
\newblock {\em Comm. Pure Appl. Math.}, 38(5):685--696, 1985.

\bibitem{W}
A.~Weinstein.
\newblock Asymptotics of eigenvalue clusters for the {L}aplacian plus a
  potential.
\newblock {\em Duke Math. J.}, 44(4):883--892, 1977.

\end{thebibliography}

\vfill\eject

\noindent D.~Bambusi,\\ 
Universit{\`a} degli studi di Milano\\
Dipartimento di Matematica\\
Via Saldini 50\\
20133 Milano, Italy\medskip\\
J.-M. Delort\\
Laboratoire Analyse G{\'e}om{\'e}trie et Applications, UMR CNRS 7539\\
Institut Galil{\'e}e, Universit{\'e} Paris-Nord,\\
99, Avenue J.-B. Cl{\'e}ment,\\
F-93430 Villetaneuse, France\medskip\\
B. Gr{\'e}bert\\
Laboratoire de Math{\'e}matiques Jean Leray, UMR CNRS 6629\\
Universit{\'e} de Nantes\\ 
2, rue de la Houssini{\`e}re\\
F-44322 Nantes Cedex 03,
France,\medskip\\
J. Szeftel\\
Department of Mathematics,\\
Princeton University,\\
Fine Hall, Washington Road\\
Princeton NJ 08544-1000 USA\\
and\\
Math{\'e}matiques Appliqu{\'e}es de Bordeaux, UMR CNRS 5466\\
Universit{\'e}  Bordeaux 1 \\
351 cours de la Lib{\'e}ration \\ 
33405 Talence cedex, France 
\end{document}